\documentclass[a4paper,leqno]{article}
\usepackage[centertags]{amsmath}
\usepackage{amsfonts}
\usepackage{amssymb}
\usepackage{amsthm}
\usepackage[active]{srcltx} 
\usepackage[utf8]{inputenc}
\newtheorem{Th}{Theorem}
\newtheorem{Main-Th}{Main Theorem}
\newtheorem{Cor}{Corollary}
\newtheorem{Pro}{Proposition}
\newtheorem{Lem}{Lemma}
\theoremstyle{definition}

\newtheorem{Exs}{Examples}

\theoremstyle{remark}
\newtheorem{Rk}{Remark}

\newenvironment{Proof}{{\bf Proof:}}
{%
\mbox{}%
\nolinebreak%
\hfill%
\rule{2mm}{2mm}%
\medbreak%
\par%
}

\numberwithin{equation}{section}

\title{\bf\large ON POINT SPECTRA OF VECTOR FIELDS}

\author{M.T.K. Abbassi and I. Lakrini}
\date{}
\begin{document}
\maketitle

\begin{abstract}
The aim of this paper is to investigate the point spectra of vector fields. We will define the point spectrum of a vector field and study some of its basic properties. In particular, we will prove that point spectra are well-behaved under the action of isometries and that the point spectra of compactly supported vector fields are trivial.



\end{abstract}


\section*{Introduction}


Vector fields are central objects in differential geometry:  they appear as sections of the tangent bundle and linear derivations of the algebra of smooth functions (cf. \cite{Ber,Lee,War}). Another reason for their importance is the fact that some classes of vector fields are infinitesimal transformations of certain geometric structures on the underlying manifold (e.g. Killing vector fields, conformal vector fields...etc). Also vector fields are infinitesimal generators of local one parameter groups of diffeomorphisms. Furthermore, the existence of vector fields which satisfies additional requirements often leads to topological or geometric obstructions on the manifold under consideration.

Vector fields are extensively studied from a geometric viewpoint. On the other hand, spectral theory has shown to be very useful in geometry and spectral techniques are now widely used in Riemannian geometry (cf. \cite{BGM,Cha,Lab}). Precisely, spectral geometry, as the study of spectral theory of the Laplace operator, is a great case of this interaction between geometry and spectral theory. It turned out that the spectrum of the Laplacian on a Riemannian manifold is a Riemannian invariant that encodes a lot of information about the geometry of that manifold. We propose here an operator theoritic approach to vector fields which will eventually lead us to consider spectral theory of vector fields.

From the viewpoint of operator theory, one of the most important spectral characteristics of a bounded linear operator, between Banach spaces, is that its spectrum, and hence the point spectrum, is bounded (see \cite{Bre,Kub,Wei,Kol-Fo} for details on bounded linear operators and their spectral theory). The present authors proved that if one makes certain choices of norms on the algebra of smooth functions on a closed Riemannian manifold, then vector fields are bounded linear operators (cf. \cite{Abb-Lak}). We have also showed that it is possible to extend these operators to appropriate Hilbert spaces in the hope of using the vast literature on spectral theory of bounded linear operators between Hilbert spaces. We shall see that the case of compact manifolds is trivial, from a spectral perspective, and that one should focus on the non-compact case. This lead us to consider the case of compactly supported vector fields on non-compact manifolds. Indeed, we extend here, in a straightforward manner, the results of the compact case to the non-compact case.

We will define the point spectrum of a vector field and we will prove that it is invariant under the action of isometries in a sense that we will precise in the sequel. More to the point, we also relate eigenfunctions of a vector field to its local flow. Moreover, we will prove that, when the manifold is compact or the vector field is compactly supported, the point spectrum is reduced to $\{0\}$. This result is extremely useful since it suggests the triviality of the compact case. Nevertheless, we shall give many examples of vector fields in the non-compact case which have non-trivial point spectra.

Unless otherwise stated, all manifolds, functions and vector fields are assumed to be smooth and smooth means differentiable of class $C^{\infty}$. We also assume that all manifolds are connected.


\section{Background material}


Let $(M,g)$ be a Riemannian manifold, a vector field $X$ on $M$ is a derivation of the algebra of smooth functions i.e. a linear mapping $X:C^{\infty}(M)\longrightarrow C^{\infty}(M)$ which satisfies the Leibniz rule. Denote by $\mathcal{D}(M)$ the space of test functions (i.e. compactly supported smooth functions).

Consider the space $L^2(M)$ of square-integrable functions on $M$. The space $L^2(M)$ is a Hilbert space with respect to the inner product
\begin{equation}\label{Inner1}
\langle f,g\rangle_{L^2}=\int_M fg \upsilon_g,
\end{equation}
for $f,g\in L^{2}(M)$. The induced norm is given by
\begin{equation}\label{NL2}
\|f\|_{L^2}=\Big(\int_M f^2 \upsilon_g\Big)^{\frac{1}{2}},
\end{equation}
for $f\in L^{2}(M)$. Further, the space $(\mathcal{D}(M),\|.\|_2)$ is dense in $L^2(M)$.

On the other hand, consider the space
\begin{equation*}
C^{2}_{1}(M)=\{f\in C^{\infty}(M):\: \|\nabla f\|\in L^{2}(M)\},
\end{equation*}
and endow this space with the norm
\begin{equation}\label{NH1}
\|f\|_{H^1}=\Big(\int_M f^2 \upsilon_g+\int_M \|\nabla f\|^2 \upsilon_g\Big)^{\frac{1}{2}},
\end{equation}
where $\upsilon_g$ is the Riemannian volume element of $(M,g)$ and $\nabla$ is its Levi-Civita connection. The Sobolev space $H^{1}(M)$ is the completion of $C^{2}_{1}(M)$ with respect to the norm $\|.\|_{H^1}$. The space $H^1(M)$ is a Hilbert space with respect to the inner product
\begin{equation}\label{Inner2}
\langle f,g\rangle_{H^1}=\int_M fg \upsilon_g+\int_M g(\nabla f,\nabla g)\upsilon_g,
\end{equation}
for $f,g\in H^1(M)$. The induced norm from this inner product is the norm $\|.\|_{H^1}$ (cf. \cite{Heb1,Heb2,Aub1} for details on the subject of Sobolev spaces on Riemannian manifolds).

The authors proved that when $(M,g)$ is a closed manifold (i.e. compact and without boundary), a vector field $X\in \mathfrak{X}(M)$, seen as a linear mapping, $X:(C^{\infty}(M),\|.\|_{H^1})\longrightarrow (C^{\infty}(M),\|.\|_{L^2})$ is a bounded linear operator with operator norm $|||X|||\leq \|X\|_\infty$, with $\|X\|_\infty=\max_{x\in M}\{\|X_x\|\}$ (cf. \cite{Abb-Lak}).

Using a density argument, one can extend the bounded linear operator $X:(C^{\infty}(M),\|.\|_{H^1})\longrightarrow (C^{\infty}(M),\|.\|_{L^2})$ to a bounded linear operator $X:H^{1}(M)\longrightarrow L^{2}(M)$ such that $|||X|||\leq \|X\|_\infty$.

In the case $M$ is non-compact, then all what is being said extend in a straightforward manner by working with compactly supported vector fields and compactly supported functions. Precisely, if $X$ be a compactly supported vector field on $M$. Since $Xf$ is a compactly supported function, for every $f\in C^{2}_{1}(M)$, the vector field $X$ induces a linear mapping $X:C^{2}_{1}(M)\longrightarrow \mathcal{D}(M)$. In fact, $X$ induces a bounded linear operator $X:C^{2}_{1}(M)\longrightarrow \mathcal{D}(M)$. Therefore, the same density argument implies the following
\begin{Pro}
Let $X$ be a compactly supported vector field, then $X$ induces a bounded linear operator $X:(H^{1}(M),\|.\|_{H^1})\longrightarrow (L^{2}(M),\|.\|_{L^2})$ such that $|||X|||\leq \|X\|_{\infty}$.
\end{Pro}


\section{Point spectra of vector fields}


Let $X\in \mathfrak{X}(M)$ be a vector field on $M$. A real number $\lambda$ is said to be an \emph{eigenvalue} of $X$ if there exists a non-zero function $f\in C^{\infty}(M)$ such that
\begin{equation}\label{Eigen-Val}
Xf=\lambda f.
\end{equation}
The set of eigenvalues of $X$ is called the \emph{point spectrum} of $X$ and we denote it by $\sigma_p(X)$. A non-vanishing function $f$ that satisfies (\ref{Eigen-Val}) for some eigenvalue $\lambda$ is called an eigenfunction of $X$. The set of eigenfunctions of $X$ associated with an eigenvalue $\lambda$ is the vector subspace of $C^{\infty}(M)$ denoted by $E_\lambda(X)$ and called the eigenspace associated with $\lambda$, it is just $\ker(X-\lambda I)$, where $I:C^{\infty}(M)\longrightarrow C^{\infty}(M)$ is the identity map.

The point spectrum of $X$ is never empty. Zero is an eigenvalue of every vector field. For the eigenspace associated to the eigenvalue 0 of a vector field $X$, we have

\begin{Pro}
   $\dim(E_0(X))$ is either 1 or $+\infty$.
\end{Pro}

\begin{Proof}
  Since constants are eigenfunctions of $X$ associated to the eigenvalue 0, then $\dim(E_0(X))\geq 1$. Suppose that $\dim(E_0(X))\geq 2$. Then there is a non-constant function $f$ which is an eigenfunction of $X$ associated to the eigenvalue 0. It is easy to see that, for each integer $k\geq1$, $f^k$ is also an eigenfunction of $X$ associated to the eigenvalue 0 ($X.f^k= f^{k-1}X.f=0$). The result follows from the following Lemma.
\end{Proof}

\begin{Lem}
Let $f \in C^\infty(M)$ be a non-constant function, then $\{f^k, k\geq 1\} $ is linearly independent.
\end{Lem}

\begin{Proof}
  Since $f$ is not constant, then there is $x_0 \in M$ such that $d_{x_0}f \neq 0$. It follows that there is $Y \in \mathfrak{X}(M)$ such that $Y_{x_0}(f) \neq 0$. By the rectification theorem of vector fields, there is a coordinate system $(U,\varphi;x^1,...,x^n)$ of $M$ around $x_0$ such that $Y\upharpoonright_{U}=\frac{\partial}{\partial x^1}$, then there is $\varepsilon >0$ such that $\frac{\partial f}{\partial x^1}(\varphi^{-1}(t,x_0^2,...,x_0^n)) \neq 0$, for all $t \in I:=]x_0^1- \varepsilon, x_0^1+ \varepsilon[$, where $(x_0^1,...,x_0^n)=\varphi(x_0)$. Consider the real function $\psi:I \rightarrow \mathbb{R}$ defined by $\psi(t):=f(\varphi^{-1}(t,x_0^2,...,x_0^n))$. The function $\psi$ is a smooth function such that $\psi^\prime(t)=\frac{\partial f}{\partial x^1}(\varphi^{-1}(t,x_0^2,...,x_0^n)) \neq 0$, for all $t\in I$. By continuity, $\psi^\prime$ maintain the same sign on $I$, and hence $\psi$ is a bijection from $I$ onto the interval $J:=\psi(I)$ of $\mathbb{R}$. Denote by $s:J \rightarrow I$ the inverse function of $\psi$.

  Now to show that $\{f^k, k \geq 1\} $ is linearly independent, consider an integer $k \geq 1$ and real numbers $\lambda_l$, $1 \leq l \leq k$ such that $\sum_{l=1}^k \lambda_l f^l=0$. Restricting the identity $\sum_{l=1}^k \lambda_l f^l=0$ to the curve image $\varphi^{-1}(I \times \{(x_0^2,...,x_0^n)\})$, we have $\sum_{l=1}^k \lambda_l \psi^l(t)=0$, for all $t \in I$, or equivalently, $\sum_{l=1}^k \lambda_l s^l=0$, for all $s \in J$, which implies that $\lambda_l=0$ for all $1 \leq l \leq k$. This completes our proof.
\end{Proof}

\begin{Exs}
\begin{itemize}$ $
\item [(i)] \textbf{Killing vector fields}: If $X$ is a Killing vector field with non-constant length, then $\dim(E_0(X))=+\infty$. Indeed, let $f$ be the smooth (non-constant) function defined by $f=\|X\|^{2}=g(X,X)$, then  $$X.f=2g(\nabla _X X,X)=0,$$
since $X$ satisfies the Killing equation (i.e. $g(\nabla_YX,Z)+g(\nabla_ZX,Y)=0$ for all $Y,Z\in \mathfrak{X}(M)$). Whence $f\in E_{0}(X)$. Using the previous proposition, we deduce that $\dim(E_0(X))=+\infty$.
\item [(ii)] \textbf{Homothetic vector fields}: A vector field $X\in \mathfrak{X}(M)$ is said to be \emph{homothetic} if its (locally defined) flow defines homothetic transformations. It is characterized by the equation $L_Xg=2cg$, for some real constant $c$. The case $c=0$ corresponds to Killing vector fields. An important example of homothetic vector fields is concurrent vector fields. A vector field $X$ is said to be a \emph{concurrent} vector field if $\nabla_YX=Y$, for all $Y\in \mathfrak{X}(M)$. In \cite{Yan}, the author proved that if the holonomy group of $(M,g)$ fixes a point, then $M$ supports a non-zero concurrent vector field (see \cite{Bri-Yan,Yan} for details). If $X$ is a non-zero homothetic vector field i.e. $L_Xg=2cg$ with $c\neq 0$ and $f$ is the function defined by $f=g(X,X)$, then
$$X.f=X.g(X,X)=2g(\nabla_XX,X)=2cg(X,X)=2cf,$$
thus $2c$ is an eigenvalue of $X$ and $f$ is an associated eigenfunction. In particular, $2$ is an eigenvalue of any non-zero concurrent vector field.

\item [(iii)]\textbf{A vector field on a non-compact submanifold of $\mathbb{R}^{2}$}: Denote by $D$ the open unit disc in $\mathbb{R}^{2}$ and let $D^*=D\setminus \{0\}$. We endow $D^*$ with the induced Riemannian metric. Consider the vector field $X\in \mathfrak{X}(D^*)$ defined by
$$X_{(x,y)}=x\frac{\partial}{\partial x}+y\frac{\partial}{\partial y},$$
for all $(x,y)\in D^*$. Let $f:D^*\longrightarrow \mathbb{R}^2$ be the function defined by $f(x,y)=\sqrt{x^2+y^2}$. Hence, for all $(x,y)\in D^*$, we have $\frac{\partial f}{\partial x}(x,y)=\dfrac{x}{f}$ and $\frac{\partial f}{\partial y}(x,y)=\dfrac{y}{f}$. We conclude that
$X.f=f$, hence $1$ is an eigenvalue and $f$ is an associated eigenfunction.
\end{itemize}
\end{Exs}

On the other hand, the boundedness of the extension of $X$ implies that the point spectrum is bounded provided that $X$ is compactly supported. As a consequence, we have the following:

\begin{Th}\label{compact-support}
  For every compactly supported vector field $X$, the point spectrum $\sigma_p(X)$ of $X$ is trivial, i.e. $\sigma_p(X)=\{0\}$. Consequently, if $M$ is compact, then $\sigma_p(X)=\{0\}$, for every $X \in \mathfrak{X}(M)$.
\end{Th}

\begin{Proof}
  Assume that there exists $\lambda \in \sigma_p(X) \setminus \{0\}$ and let $f$ be an associated eigenfunction. Then for every integer $k \geq 1$, using the Leibniz formula, we prove that $X.f^k=k \lambda f^k$, i.e. $k\lambda \in \sigma_p(X)$. We deduce that $\{k\lambda, k \geq 1\} \subset \sigma_p(X)$ which contradicts the boundedness of $\sigma_p(X)$.
\end{Proof}

\begin{Rk}
  The second part of Theorem \ref{compact-support} can be proved by a variational method as follows. If $\lambda$ is a non-vanishing eigenvalue, then there exists $f\in C^{\infty}(M)$, $f\neq 0$, such that $Xf=\lambda f$. Since $M$ is compact, $f$ is bounded and attains its bounds. In particular, there exist $x_1,x_2\in M$ such that
\begin{equation*}
f(x_1)=\max_{x\in M}\{f(x)\}\quad \mbox{and}\quad f(x_2)=\min_{x\in M}\{f(x)\}.
\end{equation*}
The point $x_1$ is critical, then
\begin{equation}
X.f(x_1)=d_{x_1}f(X_{x_1})=\lambda f(x_1)=0,
\end{equation}
which implies that $f(x_1)=0$, whence $f\leq 0$. The same argument applied to $x_2$ yields $f\geq 0$. Hence $f=0$, which is impossible.
\end{Rk}

\begin{Cor}
Assume that $(M,g)$ is compact, then $(M,g)$ does not support any non-zero concurrent vector field. Further, on a compact manifold, every non-zero homothetic vector field is a Killing vector field.
\end{Cor}

The point spectra of vector fields are well-behaved under the action of isometries. Precisely, we have:
\begin{Th}
Let $\phi:(M,g)\longrightarrow (N,h)$ is an isometry and $X\in \mathfrak{X}(M)$. Denote by $Y=\phi_*X$. Then
\begin{itemize}
\item [(i)] $\sigma_p(X)=\sigma_p(Y)$;
\item [(ii)] The isometry $\phi$ induces isomorphisms $\phi_\lambda:E_{\lambda}(X)\longrightarrow E_{\lambda}(Y)$ given by $\phi_\lambda(f)=f\circ \phi^{-1}$, for all $\lambda\in \sigma_p(X)$ and $f\in E_\lambda(X)$.
\end{itemize}
\end{Th}
\begin{Proof}
Let $\lambda\in \sigma_p(X)$ and $f\in E_\lambda(X)$, then
\begin{align*}
d(f\circ \phi^{-1})(V)&=df((\phi^{-1})_*V)\\
&=g(\nabla^g f,(\phi^{-1})_*V)\\
&=(\phi^{-1})^*g(\phi_* (\nabla^g f),V)\\
&=h(\phi_* (\nabla^g f),V),
\end{align*}
for all $V\in \mathfrak{X}(N)$. Thus $\nabla^h(f\circ \phi^{-1})=\phi_*(\nabla^g f)$. This implies that
\begin{align*}
Y.(f\circ \phi^{-1})&=h\big(Y,\nabla^h(f\circ \phi^{-1})\big)\\
&=h\big(\phi_*X,\phi_*(\nabla^g f)\big)\\
&=\phi^*h(X,\nabla^g f)\\
&=g(X,\nabla^g f)\circ \phi^{-1}\\
&=(Xf)\circ \phi^{-1},
\end{align*}
which implies at once both assertions.
\end{Proof}

Next, we explore the relationship between the spectrum and the flow of a vector field. When $M$ is compact, every vector field is complete. Thus the flow is global.

If $\gamma:(-\epsilon,\epsilon)\longrightarrow M$ is an integral curve of $X$ with initial data $\gamma(0)=x$, then
\begin{equation}\label{Int-eigen}
X.f=\overset{.}{\gamma}(t).f=\frac{d}{dt}(f\circ \gamma)(t)=\lambda f\circ \gamma(t),
\end{equation}
for every $\lambda \in \sigma_p(X)$, $f\in E_\lambda(X)$ and $t\in (-\epsilon,\epsilon)$. This implies that
\begin{equation}\label{Int-eigen2}
f\circ \gamma(t)=A\exp(\lambda t),
\end{equation}
for some real constant $A$. The initial condition implies that $A=f(x)$. Denote by $\Phi:\mathbb{R}\times M\longrightarrow M$ the local flow of $X$.
\begin{Pro}
Let $\lambda\in \sigma_p(X)$, then a function $f:M\longrightarrow \mathbb{R}$ is an eigenfunction associated with $\lambda$ if and only if
\begin{equation*}
f\circ \Phi(t,x)=f(x)\exp(\lambda t),
\end{equation*}
for every $(t,x)$ in a neighbourhood of $(0,x)$.
\end{Pro}
\begin{Proof}
Fix a point $x\in M$, since $\Phi(.,x)$ is an integral curve of $X$, then Equation (\ref{Int-eigen2}) implies $f\circ \Phi(t,x)=f(x)\exp(\lambda t)$, for every $t\in (-\epsilon,\epsilon)$ for some $\epsilon>0$. Since the point $x$ is arbitrary, the equation follows. The converse is straightforward.
\end{Proof}

\begin{Pro}
Let $X$ and $Y$ be two vector fields on $(M,g)$. If $[X,Y]=0$, then $Y(E_\lambda(X))\subseteq E_\lambda(X)$, for all $\lambda\in \sigma_p(X)$.
\end{Pro}
\begin{Proof}
Let $f\in E_\lambda(X)$, then
\begin{equation*}
XY(f)=X(Yf)=Y(Xf)=Y(\lambda f)=\lambda Yf.
\end{equation*}
\end{Proof}

In the general context, weird phenomena can occur.  We give here a striking example of a non-compactly supported vector field on a non-complete Riemannian manifold for which the point spectrum $\sigma_p(X)=\mathbb{R}$. Denote by $$L=\{(0,y):\: y\in \mathbb{R}\},$$
and consider the open submanifold $\Omega=\mathbb{R}^2\setminus L$ endowed with the induced Riemannian metric by the usual Euclidean metric on the plane.

Define the vector field $X \in \mathfrak{X}(\Omega)$ by setting $$X_{(x,y)}=-y\frac{\partial }{\partial x}+x\frac{\partial }{\partial y},$$
for all $(x,y)\in \Omega$.

For $\lambda\in \mathbb{R}^*$, consider the function $\varphi_\lambda:\Omega \longrightarrow \mathbb{R}$ given by $$\varphi_\lambda(x,y)=\exp\big(\lambda \arctan(\frac{y}{x})\big),$$ which is smooth on $\Omega$. Hence
\begin{align*}
\frac{\partial \varphi_\lambda}{\partial x}(x,y)&=-\frac{\lambda y}{x^2+y^2}\varphi_\lambda(x,y),\\
\frac{\partial \varphi_\lambda}{\partial y}(x,y)&=\frac{\lambda x}{x^2+y^2}\varphi_\lambda(x,y),
\end{align*}
whence
\begin{equation*}
X.\varphi_\lambda=\lambda \varphi_\lambda.
\end{equation*}
Thus $\lambda\in \sigma_p(X)$ and $\varphi_\lambda\in E_\lambda(X)$. Thus $\sigma_p(X)=\mathbb{R}$.
\begin{Rk}
Let $(M,g)$ a Riemannian manifold and $X \in \mathfrak{X}(M)$. For $\lambda\in \sigma_p(X)$, the equation $Xf=\lambda f$ has a local character. More precisely, let $(U,x^1,...,x^n)$ is a local coordinate system of $M$. If the vector field $X$ is expressed in this local coordinate system as $X=\sum_{i=1}^{n}X^i\frac{\partial}{\partial x^i}$, then $Xf=\lambda f$ becomes $\sum_{i}X^i\frac{\partial f}{\partial x^i}=\lambda f$. Conversely, if the equation holds in every local coordinate system, then it holds globally.
\end{Rk}


\end{document}